\catcode`\@=11

\font\tenmsa=msam10 \font\sevenmsa=msam7 \font\fivemsa=msam5
\font\tenmsb=msbm10
\font\sevenmsb=msbm7 \font\fivemsb=msbm5 \newfam\msafam \newfam\msbfam
\textfont\msafam=\tenmsa \scriptfont\msafam=\sevenmsa
\scriptscriptfont\msafam=\fivemsa \textfont\msbfam=\tenmsb
\scriptfont\msbfam=\sevenmsb \scriptscriptfont\msbfam=\fivemsb

\def\hexnumber@#1{\ifnum#1<10 \number#1\else \ifnum#1=10 A\else\ifnum#1=11
 B\else\ifnum#1=12 C\else \ifnum#1=13 D\else\ifnum#1=14 E\else\ifnum#1=15
 F\fi\fi\fi\fi\fi\fi\fi}

\def\msa@{\hexnumber@\msafam} \def\msb@{\hexnumber@\msbfam}
\mathchardef\boxdot="2\msa@00 \mathchardef\boxplus="2\msa@01
\mathchardef\boxtimes="2\msa@02 \mathchardef\square="0\msa@03
\mathchardef\blacksquare="0\msa@04 \mathchardef\centerdot="2\msa@05
\mathchardef\lozenge="0\msa@06 \mathchardef\blacklozenge="0\msa@07
\mathchardef\circlearrowright="3\msa@08 \mathchardef\circlearrowleft="3\msa@09
\mathchardef\rightleftharpoons="3\msa@0A
\mathchardef\leftrightharpoons="3\msa@0B \mathchardef\boxminus="2\msa@0C
\mathchardef\Vdash="3\msa@0D \mathchardef\Vvdash="3\msa@0E
\mathchardef\vDash="3\msa@0F \mathchardef\twoheadrightarrow="3\msa@10
\mathchardef\twoheadleftarrow="3\msa@11 \mathchardef\leftleftarrows="3\msa@12
\mathchardef\rightrightarrows="3\msa@13 \mathchardef\upuparrows="3\msa@14
\mathchardef\downdownarrows="3\msa@15 \mathchardef\upharpoonright="3\msa@16
 \mathchardef\downharpoonright="3\msa@17
\mathchardef\upharpoonleft="3\msa@18 \mathchardef\downharpoonleft="3\msa@19
\mathchardef\rightarrowtail="3\msa@1A \mathchardef\leftarrowtail="3\msa@1B
\mathchardef\leftrightarrows="3\msa@1C \mathchardef\rightleftarrows="3\msa@1D
\mathchardef\Lsh="3\msa@1E \mathchardef\Rsh="3\msa@1F
\mathchardef\rightsquigarrow="3\msa@20
\mathchardef\leftrightsquigarrow="3\msa@21 \mathchardef\looparrowleft="3\msa@22
\mathchardef\looparrowright="3\msa@23 \mathchardef\circeq="3\msa@24
\mathchardef\succsim="3\msa@25 \mathchardef\gtrsim="3\msa@26
\mathchardef\gtrapprox="3\msa@27 \mathchardef\multimap="3\msa@28
\mathchardef\therefore="3\msa@29 \mathchardef\because="3\msa@2A
\mathchardef\doteqdot="3\msa@2B 
\mathchardef\traceiangleq="3\msa@2C \mathchardef\precsim="3\msa@2D
\mathchardef\lesssim="3\msa@2E \mathchardef\lessapprox="3\msa@2F
\mathchardef\eqslantless="3\msa@30 \mathchardef\eqslantgtr="3\msa@31
\mathchardef\curlyeqprec="3\msa@32 \mathchardef\curlyeqsucc="3\msa@33
\mathchardef\preccurlyeq="3\msa@34 \mathchardef\leqq="3\msa@35
\mathchardef\leqslant="3\msa@36 \mathchardef\lessgtr="3\msa@37
\mathchardef\backprime="0\msa@38 \mathchardef\risingdotseq="3\msa@3A
\mathchardef\fallingdotseq="3\msa@3B \mathchardef\succcurlyeq="3\msa@3C
\mathchardef\geqq="3\msa@3D \mathchardef\geqslant="3\msa@3E
\mathchardef\gtrless="3\msa@3F \mathchardef\sqsubset="3\msa@40
\mathchardef\sqsupset="3\msa@41
\mathchardef\trianglelefteq="3\msa@45 \mathchardef\bigstar="0\msa@46
\mathchardef\between="3\msa@47 \mathchardef\blacktriangledown="0\msa@48
\mathchardef\blacktriangleright="3\msa@49
\mathchardef\blacktriangleleft="3\msa@4A
\mathchardef\blacktriangle="0\msa@4E \mathchardef\triangledown="0\msa@4F
\mathchardef\eqcirc="3\msa@50 \mathchardef\lesseqgtr="3\msa@51
\mathchardef\gtreqless="3\msa@52 \mathchardef\lesseqqgtr="3\msa@53
\mathchardef\gtreqqless="3\msa@54 \mathchardef\Rrightarrow="3\msa@56
\mathchardef\Lleftarrow="3\msa@57 \mathchardef\veebar="2\msa@59
\mathchardef\barwedge="2\msa@5A \mathchardef\doublebarwedge="2\msa@5B
\mathchardef\angle="0\msa@5C \mathchardef\measuredangle="0\msa@5D
\mathchardef\sphericalangle="0\msa@5E \mathchardef\varpropto="3\msa@5F
\mathchardef\smallsmile="3\msa@60 \mathchardef\smallfrown="3\msa@61
\mathchardef\Subset="3\msa@62 \mathchardef\Supset="3\msa@63
\mathchardef\Cup="2\msa@64  \mathchardef\Cap="2\msa@65
 \mathchardef\curlywedge="2\msa@66
\mathchardef\curlyvee="2\msa@67 \mathchardef\leftthreetimes="2\msa@68
\mathchardef\rightthreetimes="2\msa@69 \mathchardef\subseteqq="3\msa@6A
\mathchardef\supseteqq="3\msa@6B \mathchardef\bumpeq="3\msa@6C
\mathchardef\Bumpeq="3\msa@6D \mathchardef\lll="3\msa@6E 
\mathchardef\ggg="3\msa@6F  \mathchardef\circledS="0\msa@73
\mathchardef\pitchfork="3\msa@74 \mathchardef\dotplus="2\msa@75
\mathchardef\backsim="3\msa@76 \mathchardef\backsimeq="3\msa@77
\mathchardef\complement="0\msa@7B \mathchardef\intercal="2\msa@7C
\mathchardef\circledcirc="2\msa@7D \mathchardef\circledast="2\msa@7E
\mathchardef\circleddash="2\msa@7F \def\ulcorner{\delimiter"4\msa@70\msa@70 }
\def\urcorner{\delimiter"5\msa@71\msa@71 }
\def\llcorner{\delimiter"4\msa@78\msa@78 }
\def\lrcorner{\delimiter"5\msa@79\msa@79 } \def\yen{\mathhexbox\msa@55 }
\def\checkmark{\mathhexbox\msa@58 } \def\circledR{\mathhexbox\msa@72 }
\def\maltese{\mathhexbox\msa@7A } \mathchardef\lvertneqq="3\msb@00
\mathchardef\gvertneqq="3\msb@01 \mathchardef\nleq="3\msb@02
\mathchardef\ngeq="3\msb@03 \mathchardef\nless="3\msb@04
\mathchardef\ngtr="3\msb@05 \mathchardef\nprec="3\msb@06
\mathchardef\nsucc="3\msb@07 \mathchardef\lneqq="3\msb@08
\mathchardef\gneqq="3\msb@09 \mathchardef\nleqslant="3\msb@0A
\mathchardef\ngeqslant="3\msb@0B \mathchardef\lneq="3\msb@0C
\mathchardef\gneq="3\msb@0D \mathchardef\npreceq="3\msb@0E
\mathchardef\nsucceq="3\msb@0F \mathchardef\precnsim="3\msb@10
\mathchardef\succnsim="3\msb@11 \mathchardef\lnsim="3\msb@12
\mathchardef\gnsim="3\msb@13 \mathchardef\nleqq="3\msb@14
\mathchardef\ngeqq="3\msb@15 \mathchardef\precneqq="3\msb@16
\mathchardef\succneqq="3\msb@17 \mathchardef\precnapprox="3\msb@18
\mathchardef\succnapprox="3\msb@19 \mathchardef\lnapprox="3\msb@1A
\mathchardef\gnapprox="3\msb@1B \mathchardef\nsim="3\msb@1C
\mathchardef\napprox="3\msb@1D
\mathchardef\nsupseteqq="3\msb@23 \mathchardef\subsetneqq="3\msb@24
\mathchardef\supsetneqq="3\msb@25
\mathchardef\supsetneq="3\msb@29 \mathchardef\nsubseteq="3\msb@2A
\mathchardef\nsupseteq="3\msb@2B \mathchardef\nparallel="3\msb@2C
\mathchardef\nmid="3\msb@2D \mathchardef\nshortmid="3\msb@2E
\mathchardef\nshortparallel="3\msb@2F \mathchardef\nvdash="3\msb@30
\mathchardef\nVdash="3\msb@31 \mathchardef\nvDash="3\msb@32
\mathchardef\nVDash="3\msb@33 \mathchardef\ntrianglerighteq="3\msb@34
\mathchardef\ntrianglelefteq="3\msb@35 \mathchardef\ntriangleleft="3\msb@36
\mathchardef\ntriangleright="3\msb@37 \mathchardef\nleftarrow="3\msb@38
\mathchardef\nrightarrow="3\msb@39 \mathchardef\nLeftarrow="3\msb@3A
\mathchardef\nRightarrow="3\msb@3B \mathchardef\nLeftrightarrow="3\msb@3C
\mathchardef\nleftrightarrow="3\msb@3D \mathchardef\divideontimes="2\msb@3E
\mathchardef\varnothing="0\msb@3F \mathchardef\nexists="0\msb@40
\mathchardef\mho="0\msb@66 \mathchardef\thorn="0\msb@67
\mathchardef\beth="0\msb@69 \mathchardef\gimel="0\msb@6A
\mathchardef\daleth="0\msb@6B \mathchardef\lessdot="3\msb@6C
\mathchardef\gtrdot="3\msb@6D \mathchardef\ltimes="2\msb@6E
\mathchardef\rtimes="2\msb@6F \mathchardef\shortmid="3\msb@70
\mathchardef\shortparallel="3\msb@71 \mathchardef\smallsetminus="2\msb@72
\mathchardef\thicksim="3\msb@73 \mathchardef\thickapprox="3\msb@74
\mathchardef\approxeq="3\msb@75 \mathchardef\succapprox="3\msb@76
\mathchardef\precapprox="3\msb@77 \mathchardef\curvearrowleft="3\msb@78
\mathchardef\curvearrowright="3\msb@79 \mathchardef\digamma="0\msb@7A
\mathchardef\varkappa="0\msb@7B \mathchardef\hslash="0\msb@7D
\mathchardef\hbar="0\msb@7E \mathchardef\backepsilon="3\msb@7F
\def\Bbb{\ifmmode\let\next\Bbb@\else
\def\next{\errmessage{Use \string\Bbb\space only in math mode}}\fi\next}
\def\Bbb@#1{{\Bbb@@{#1}}} \def\Bbb@@#1{\fam\msbfam#1}
\catcode`\@=\active

\def\CR{\hbox{{$\cal R$}}}

\def\cu{\hbox{\sl u}} % used for special element u
\def\cg{\hbox{{\sl g}}} % used for Lie algebra 'gothic g'
\def\lform{\hbox{$\sqcup$}\llap{\hbox{$\sqcap$}}}
\def\h{{{1\over2}}}

\def\R{{\Bbb R}}
\def\C{{\Bbb C}}
\def\Z{{\Bbb Z}}

\def\eps{{\epsilon}}

\def\dcross{{\bowtie}}

\def\rbiprod{{\cdot\kern-.33em\triangleright\!\!\!<}}
\def\lbiprod{{>\!\!\!\triangleleft\kern-.33em\cdot\, }}

\def\tens{\mathop{\otimes}}

\def\la{{\triangleright}}\def\ra{{\triangleleft}}

\def\isom{{\cong}}

\def\ev{{\rm ev}}

\def\id{{\rm id}}

\def\<{\langle}
\def\>{\rangle}
\def\dila{{\varsigma}}

\def\equad{\kern -1.7em}

\def\eqn#1#2{\begin{equation}#2\label{#1}\end{equation}}

\def\o{{}_{\scriptscriptstyle(1)}}
\def\t{{}_{\scriptscriptstyle(2)}}
\def\th{{}_{\scriptscriptstyle(3)}}

\def\bo{{}^{\bar{\scriptscriptstyle(1)}}}
\def\bt{{}^{\bar{\scriptscriptstyle(2)}}}
\def\Ro{{\CR^{\scriptscriptstyle(1)}}}
\def\Rt{{\CR^{\scriptscriptstyle(2)}}}

\def\und#1{{\underline {#1}}}

\def\uo{{{}^{\scriptscriptstyle(1)}}}
\def\ut{{{}^{\scriptscriptstyle(2)}}}

\def\umo{{{}^{\scriptscriptstyle-(1)}}}
\def\umt{{{}^{\scriptscriptstyle-(2)}}}

\def\Bo{{{}_{\und{\scriptscriptstyle(1)}}}}
\def\Bt{{{}_{\und{\scriptscriptstyle(2)}}}}
\def\Bth{{{}_{\und{\scriptscriptstyle(3)}}}}

\def\text#1{\mbox{\rm #1}}
\def\note#1{}

\def\blacksquare{{\lform}}%AMS Tex Fakes
\def\frac#1#2{{{#1\over#2}}}

\def\proof{\goodbreak\noindent{\bf Proof\quad}}

\def\endproof{{\ $\lform$}\bigskip }

\def\align#1{\begin{eqnarray*}#1\end{eqnarray*}}
\def\alignn#1#2{\begin{eqnarray}\label{#1}#2
\end{eqnarray}}
\def\cmath#1{\[\begin{array}{c} #1 \end{array}\]}
\def\ceqn#1#2{\begin{equation}\label{#1}\begin{array}{c}#2
\end{array}\end{equation}}

\def\vecm{{\bf m}}\def\vece{{\bf e}}
\def\vecf{{\bf f}}

\def\Ro#1{{\CR_{#1}\uo}}
\def\Rt#1{{\CR_{#1}\ut}}
\def\Rmo#1{{\CR_{#1}\umo}}
\def\Rmt#1{{\CR_{#1}\umt}}

\documentstyle{article}
\newtheorem{lemma}{Lemma}[section] \newtheorem{propos}[lemma]{Proposition}
 \newtheorem{theorem}[lemma]{Theorem}

\begin{document}\baselineskip 15pt

{\ }\qquad\qquad \hskip 1.3in  Submitted Proc. Quantum Groups and Physics, 
Prague, June 1996
\vspace{.2in}

%%%%%%%%%%% here starts the paper itself %%%%%%%%%%%%%%%%%%%%%

\begin{center} {\LARGE NEW QUANTUM GROUPS BY DOUBLE-BOSONISATION}
\\ \baselineskip 13pt{\ }
{\ }\\
S. Majid\footnote{Royal Society University Research Fellow and Fellow of
Pembroke College, Cambridge}\\
{\ }\\
Department of Mathematics, Harvard University\\
Science Center, Cambridge MA 02138, USA\footnote{During the calendar years 1995
+ 1996}\\
+\\
Department of Applied Mathematics \& Theoretical Physics\\
University of Cambridge, Cambridge CB3 9EW\\
\end{center}

\vspace{20pt}
\begin{quote}\baselineskip 11pt
\noindent{\bf Abstract} We obtain new family of quasitriangular Hopf
algebras $\C^{0|n}_q\lbiprod \widetilde{U_q(su_n)}\rbiprod \C^{0|n}_q$
via the author's recent double-bosonisation construction for new
quantum groups. They are versions of $U_q(su_{n+1})$ with a fermionic
rather than bosonic quantum plane of roots adjoined to $U_q(su_n)$. We
give the $n=2$ case in detail.  We also consider the anyonic-double of
an anyonic ($\Z/n\Z$-graded) braided group and the double-bosonisation
of the free braided group in $n$ variables.


\end{quote}
\baselineskip 13.1pt

\section{Introduction}

There are two well-known general constructions for quasitriangular Hopf
algebras or `strict quantum groups' (with universal R-matrix), namely (i)
Drinfeld's quantum double\cite{Dri} of any Hopf algebra and (ii) the R-matrix
approach introduced for general R-matrices  in \cite{Ma:mor}, where
the existence of a universal R-matrix functional on $A(R)$ was established for
for the first time (this is not a result to be found in the standard FRT
work\cite{FRT:lie}; it is needed when the universal R-matrix is not known by
other means). Recently, in \cite{Ma:dbos}, we introduced a third and
more powerful general construction for quasitriangular Hopf algebras called
{\em double-bosonisation}. It associates to any braided group $B$ covariant
under a background quantum group $H$ a quasitriangular Hopf algebra
$B^{\star\rm op}\lbiprod H\rbiprod B$ built on $B^{\star}\tens H\tens B$ and
consisting of $H$ extended by $B$ as additional `positive roots' and its dual
$B^{\star}$ as additional `negative roots'. The construction is more powerful
than the quantum double because one can reach directly to the usual $U_q(\cg)$
(the quantum double is too big and one has to quotient it).
As an application, we can now build up $U_q(\cg)$ by induction, adjoining to
$U_q(su_2)$ a quantum-braided plane of roots to get to $U_q(su_3)$ etc. See
\cite{Ma:dbos}. One can also go up the $U_q(so_n)$ series etc., for example
constructing the conformal algebra $U_q(so_6)$ from $U_q(so_4)$ by adjoining a
q-Euclidean braided plane of roots\cite{Ma:conf}.

Moreover, because each double-bosonisation has a decomposition into its three
tensor factors, we obtain an `inductive block-basis' for each family of the
classical families of quantum groups. For the $A_n$ series it is:
\eqn{An}{ \cdots \C_q^n{}^{\star\rm op}\lbiprod
(\C_q^{(n-1)}{}^{\star\rm op}\lbiprod ...\lbiprod(\C_q{}^{\star\rm
op}\lbiprod U(1)\rbiprod \C_q)U(1)\rbiprod
...U(1)\rbiprod \C_q^{n-1})U(1)\rbiprod \C_q^n\cdots}
where we can go up to arbitrary $n$. This shows $U_q(su_{n+1})$ built up by
starting with $U(1)$, adjoinig a braided line $\C_q$ to obtain
$U_q(su_2)=\C_q{}^{\star\rm op}\lbiprod U(1)\rbiprod \C_q$, then adjoining a
braided plane $\C_q^2$ to obtain $U_q(su_3)$ from this, etc. Every element of
$U_q(su_{n+1})$ decomposes into the product of a unique element of the tensor
product of these blocks.  Also, wherever we take the brackets (on either side)
we obtain a Hopf algebra. For example,
\eqn{Anpara}{ (\C_q^{(n-1)}{}^{\star\rm op}\lbiprod
...\lbiprod(\C_q{}^{\star\rm op}\lbiprod U(1)\rbiprod \C_q)U(1)\rbiprod
...U(1)\rbiprod
\C_q^{n-1})U(1)\rbiprod \C_q^n=\widetilde{U_q(su_n)}\rbiprod \C_q^n}
is a sub-Hopf algebra of $U_q(su_{n+1})$, its natural `maximal parabolic'
inhomogeneous quantum group. Thus we obtain not only the classical
q-deformations but a natural description of how they are built up from each
other. Thus, double-bosonisation is the quantum analogue of adjoining a node to
a Dynkin diagram and thereby extending its Lie algebra to a larger one.
Moreover, fixing bases of each of the $\C_q^n$, we obtain a concrete {\em
inductive basis} for the entire $A$-series quantum groups. Similarly for the
other series of $U_q(\cg)$.

On the other hand, it is clear that this construction can be used just as well
to obtain non-standard quantum groups. Basically, the double-bosonisation
construction generates not a line but a {\em tree of quantum groups}: At each
node of the tree we have the choice to adjoin any braided group covariant under
the quantum group at that node. Some of the nodes in the tree will be classical
in the sense that they have smooth $q\to 1$ limits and some others will be
purely quantum with no limit. For example, at the node $U_q(su_2)$ we can
choose $\C_q^2$ giving $U_q(su_3)$ as above, $\R_q^3$ giving $U_q(so_5)$,
$\C^{0|2}_q$ giving a non-standard quantum group without a classical limit, and
probably other choices as well. From each of these nodes we have still more
choices.  The enumeration of all braided groups covariant under each the
quantum group at each node leads to an enlarged `Dynkin diagram system'
describing the tree generated by all possible double-bosonisations.

In Section~3 we demonstrate this theory with an explicit computation of the
non-standard quantum group obtained from the node $U_q(su_2)$ by adjoining
$\C^{0|2}_q$. There is obviously an infinite number of new quantum groups to be
obtained this way, and a combinatorial challenge to elaborate the full tree
structure. In particular, at each $U_q(su_n)$ node we have a fermionic partner
$\C_q^{0|n}$ to the bosonic quantum plane $\C_q^n$, giving a fermionic cousin
of $U_q(su_{n+1})$. Moreover, being an abstract construction, double
bosonisation is not tied to any $q$ at all, and works as well for discrete
quantum groups or at roots of unity.
We give some examples of this type in Section~4.

The paper begins in Section~2, where we recall the abstract bosonisation
construction from \cite{Ma:dbos}. Full proofs are in \cite{Ma:dbos}, i.e. here
we announce the main results only. In fact,  \cite{Ma:dbos} was rejected after
several months with Journal of Algebra  on rather trivial grounds, requiring
resubmission elsewhere.

\section{Abstract Double-Bosonisation Construction}

This preliminaries section announces/recalls the abstract double-bosonisation
construction from \cite{Ma:dbos}, as needed for the construction of the two
examples in later sections. The reader can also read the latter first and keep
this section for reference as a brief introduction to \cite{Ma:dbos}.

We assume in this section that the reader is familiar with the basic ideas of
quantum groups and braided groups. We use the notations from the textbook
\cite{Ma:book}. Recall only that a quantum group $(H,\CR)$ generates a braiding
$\Psi:B\tens B\to B\tens B$, and a braided group $B$ means an $H$-covariant
algebra and coalgebra which  are compatible in the braided
sense\cite{Ma:introp} that the coproduct $\und\Delta:B\to B\und\tens B$ is an
algebra homomorphism. Here $B\und\tens B$ has the braided tensor product
algebra structure
\eqn{btens}{ (b\tens d)(c\tens e)=b\Psi(d\tens c)e,\quad \Psi(d\tens
c)=c\ra\CR\uo\tens d\ra\CR\ut,\quad \forall b,c,d,e\in B.}
Let $B^\star$ denote a braided group dual to $B$ in the braided sense of an
evaluation pairing $\ev:B^\star\tens B\to \C$ being given\cite{Ma:introp}.
This section works over a general field in place of $\C$ here, or (with
suitable care) over a commutative ring such as $\C[[\hbar]]$. We take
conventions with all quantum group actions $\ra$ from the right. The main
formulae with left actions $\la$ are summarized in the appendix of
\cite{Ma:conf}.

\begin{theorem}\cite{Ma:dbos} There is a unique Hopf algebra structure
$B^{\star\rm op}\lbiprod H\rbiprod B$ on
$B^\star\tens H\tens B$, the {\em double-bosonisation}, containing $B^{\star\rm
op},H,B$ as subalgebras, with the cross relations and coproduct
\cmath{ bh=h\o (b\ra h\t),\quad ch=h\t (c\ra h\o),\quad
 b\Bo \Ro{} c\Bo \ev(c\Bt\ra \Rt{}, b\Bt)=\ev(c\Bo,b\Bo\ra\Ro{}) c\Bt\Rt{}
b\Bt\\
\Delta b=b\Bo\ra \Ro{}\tens \Rt{} b\Bt,\quad \Delta c=\Ro{}c\Bo
\tens c\Bt\ra\Rt{}}
for all $b\in B$, $c\in B^\star$, $h\in H$.
\end{theorem}

The paper \cite{Ma:dbos} provides several other conventions in which this can
be presented (in particular, the above is {\em not} the most natural from the
point of view of a left-right symmetry interchanging the roles of $B^\star,B$).
It also proves several results about the double-bosonisation. The first is that
$H,B$ and $B^{\star\rm op},H$ generate sub-Hopf algebras $H\rbiprod B$ and
$B^{\star\rm op}\lbiprod H$ respectively. These are usual {\em bosonisations}
associated to any braided group by the author's earlier bosonisation
construction\cite{Ma:bos}. The second is that,
being built on the triple tensor product, there is an explicit product law
\alignn{dbosprod}{ (c\tens h\tens b)\cdot(d\tens g\tens a)\equad&&=
(d\Bt\ra\Rmo3 Sh\o)c\tens h\t \Ro1 \Rmt4 g\o\tens (b\Bt\ra\Ro2 g\t)a
\nonumber\\
&&\quad\ev( d\Bo\ra\Rmo4\Rmo5,b\Bo)\ev(  \und S^{-1}
d\Bth\ra\Rmt3\Rt1\Rmt5\Rt2,b\Bth)}
between general elements $c,d\in B^\star$, $h,g\in H$ and $b,a\in B$. Here
$\CR_1,\cdots,\CR_5$ are five copies of the universal R-matrix of $H$. The
third result is that when $B$ has a basis $\{e_a\}$ with dual basis $\{f^a\}$
of $B^\star$ (a strict duality) then the double-bosonisation has a universal
R-matrix
\ceqn{dbos-R}{ \CR=\exp_B{}^{-1}_{21}\cdot \CR,\quad
\exp_B{}^{-1}_{21}=f^a\tens \und S e_a.}
The fourth is that the double-bosonisation acts covariantly on $B$ itself. The
(right) action is
\eqn{dbos-act}{ v\ra b=(\und S b\Bo\ra\Ro{})(v\ra\Rt{})b\Bt,\quad v\ra
c=\ev(\und S^{-1}
c,v\Bt\ra\Rmt{})v\Bo\ra\Rmo{} }
for all $v,b\in B$ and $c\in B^\star$, and the given action $\ra$ of $H$ on
$B$.

The paper \cite{Ma:dbos} also relates the double-bosonisation to the quantum
double of
the usual bosonisation. We adopt conventions in which $D(H)=H^{*\rm op}\dcross
H$. More precisely, we suppose that $H^*$ is a Hopf algebra dually paired with
$H$ and we use the generalised quantum double\cite{Ma:mor} in the
infinite-dimensional or degenerately paired case. We also assume that the
action of $H$ on $B$ is given by evaluation against a  coaction  of $H^*$. The
paper\cite{Ma:dbos} gives the precise formula for $D(B\lbiprod H)$ where $B$ is
left-covariant under $H$. In our case, $B^\star$ is right-covariant under $H$
and hence left-covariant under $H^{\rm op}$, and we need $D(B^\star\lbiprod
H^{\rm op})$. This contains $B^{\star\rm op},H,H^*,B$ as subalgebras, with the
cross relations and coproduct\cite{Ma:dbos}
\ceqn{D(bos)}{ ah=h\t a\t\<h\th,a\o\>\<h\o,S^{-1}a\th\>,\quad c h=h\t(c\ra
h\o),\quad  b a =a\t (b\ra\Rt{}) \<\Ro{},a\o\>\\ a c =(c\ra\Rt{}) a\o
\<\Ro{},S^{-1}a\t\>,\quad  b h = h\o (b\ra h\t)\\
 b c  =\Ro{1}c\Bt  (b\Bo\ra\Ro{2})\bo (b\Bt\ra\Ro{3})
\ev(c\Bo,(b\Bo\ra\Ro{2})\bt)\ev(\und
S^{-1}c\Bth\ra\Rt{1}\Rt{2}\Rt{3},b\Bth)\\
\Delta c =\Ro{}c\Bo  \tens c\Bt\ra\Rt{},\quad \Delta b=
b\Bo\bt\tens  b\Bo\bo b\Bt}
for all $c\in B^\star$, $b\in B$, $h\in H$ and $a\in H^*$, and the usual
coproducts of $H,H^{*\rm cop}$. Here the quantum double of $H$ in the form
$H\dcross H^{*\rm
cop}$ appears as a sub-Hopf algebra. In the case of strict duality we have
Drinfeld's quasitriangular structure as
\eqn{D(bos)-R}{ \CR_D=(\und S f^a\ra (S\Rt{})\cu f^\alpha\tens \Ro{}
f^\beta\tens 1\tens
1)\tens(1\tens 1\tens e_\alpha e_\beta\tens e_a)}
where $\{e_a\}$ is a basis of $B$ with dual $\{f^a\}$ and $\{e_\alpha\}$ a
basis of $H^*$ with dual $\{f^\alpha\}$.

\begin{theorem}\cite{Ma:dbos} There is a Hopf algebra surjection
$\pi:D(B^\star\lbiprod H^{\rm op}) \to B^{\star\rm op}\lbiprod H\rbiprod B$
with
\[\pi(c)=c,\quad \pi(h)=h,\quad \pi(a)= \Rt{}\<a,\Ro{}\>,\quad
\pi(b)=b,\quad\forall c\in B^\star,\ h\in H,\ a\in H^*,\ b\in B. \]
In the case of strict duality this is a surjection of
quasitriangular Hopf algebras.
\end{theorem}

Finally, it is known\cite{Ma:skl} that bosonisations can be viewed as examples
of a more general  biproduct construction (they are not the same and did not
arise this way, however; see\cite{Ma:com}); likewise, given in \cite{Ma:dbos}
is a more general `double-biproduct' construction. The input data is a braided
group $B$  right-covariant under the quantum double $D(H)$, where $H$ is a
general Hopf algebra. It is natural to demand that $H$ has a skew-antipode but
this is not actually needed (one does not need the inverse of the braiding for
the following results). Also, it is standard to write covariance under $D(H)$
slightly more generally as a (right) crossed $H$-modules in the standard way;
see \cite{Ma:book}. This means a compatible right action and coaction of $H$ on
$B$, under which it is covariant and forms a braided group. We let $B^\odot$ be
a braided group which is {\em left}-covariant under $D(H)$ (a left crossed
$H$-module)  and `skew-dual' to $B$ in the sense\cite{Ma:dbos}
\ceqn{bipair}{ b\bo\ra c\bo\tens b\bt\la c\bt=b\tens c,\quad \sigma(h\la
c,b)=\sigma(c,b\ra h)\\
\sigma(c,ab)=\sigma(c\Bt,a\ra b\bt)\sigma(c\Bo,b\bo),\quad
\sigma(cd,b)=\sigma(c,b\Bo)\sigma(d,b\Bt),\quad \sigma(\und S
c,b)=\sigma(c,\und S^{-1}b)}
for all $h\in H$, $a,b\in B$ and $c,d\in B^\odot$, for some linear map
$\sigma:B^\odot\tens B\to \C$.

\begin{theorem}\cite{Ma:dbos} For $B^\odot,H,B$ as described, there is a unique
Hopf algebra $B^\odot\lbiprod H\rbiprod B$ built on $B^\odot\tens H\tens B$,
the
{\em double biproduct}, containing $B^\odot,H,B$ as subalgebras, with the cross
relations and coproduct
\cmath{ bh=h\o(b\ra h\t),\quad hc=(h\o\la c)h\t,\quad
 bc=b\Bo\bt c\Bt b\Bt c\Bth\bo \sigma(c\Bo,b\Bo\bo)\sigma(\und S
c\Bth\bt,b\Bth\>\\
 \Delta b=b\Bo\bo\tens b\Bo\bt b\Bt,\quad \Delta c=c\Bo c\Bt\bo\tens c\Bt\bt}
for $b\in B$, $c\in B^\odot$ and $h\in H$, and the usual coproduct of $H$.
\end{theorem}
This reduces to Theorem~2.1 in the case where $H$ is quasitriangular; and
reduces to a third construction\cite{Ma:dbos} when $H$ is dual-triangular.

\section{A Fermionic Cousin of $U_q(su_3)$}

In this section we give a non-standard example of the double-bosonisation
construction. That is, we use the R-matrix version of the double-bosonisation
theorem
in \cite{Ma:dbos} and compute it for a non-standard R-matrix. This R-matrix
double-bosonisation is built from a quantum group $H$ with generators
$\vecm^\pm=\{m^\pm{}^i{}_j\}, \dila$ obeying
\ceqn{matH}{R\vecm^\pm_1\vecm^\pm_2=\vecm^\pm_2\vecm^\pm_1R,\quad
R\vecm^+_1\vecm^-_2=\vecm^-_2\vecm^+_1R,\quad
\Delta\vecm^\pm=\vecm^\pm\tens\vecm^\pm,\quad \eps\vecm^\pm=\id\\
{}[\dila,\vecm^\pm]=0,\quad \Delta\dila=\dila\tens\dila,\quad\eps\dila=1}
and braided groups $B,B^{\star}$ with generators $\vece=\{e^i\}$ and
$\vecf=\{f_i\}$ respectively, obeying
\ceqn{matBC}{\vece_2\vece_1=R'\vece_1\vece_2,\quad
\Psi(\vece_2\tens\vece_1)=R\vece_1\tens\vece_2,\quad \und\Delta\vece=\vece\tens
1+1\tens\vece,\\
\vecf_2\vecf_1=\vecf_1\vecf_2 R',\quad \Psi(\vecf_2\tens\vecf_1)=\vecf_1\tens
\vecf_2 R,\quad\und\Delta\vecf=\vecf\tens 1+1\tens\vecf}
This is a right handed setting of the author's R-matrix braided group theory;
see \cite[Chapter~10]{Ma:book} for more details and history. There are further
relations beyond (\ref{matH}) to form a Hopf algebra, and the choice of an
associated {\em quantum group normalisation constant} $\lambda$ such that this
is quasitriangular. The double-bosonisation in this case
has cross relations and coproduct\cite{Ma:dbos}
\ceqn{matdbos}{\vece_2\vecm^+_1=\lambda R\vecm^+_1\vece_2,\quad
\vecm^-_2\vece_1=\lambda R\vece_1\vecm^-_2,\quad
\vecm^+_1\vecf_2=\vecf_2\vecm^+_1 \lambda R,\quad
\vecf_1\vecm^-_2=\vecm^-_2\vecf_1 \lambda R\\
 \dila\vecf=\lambda\vecf\dila,\quad \vece\dila=\lambda\dila\vece,\quad
[\vece,\vecf]={\vecm^+\dila^{-1}-\dila\vecm^-\over q-q^{-1}}\\
\Delta e^i=e^a\tens m^+{}^i{}_a\dila^{-1}+1\tens e^i,\quad \Delta f_i=f_i\tens
1+\dila m^-{}^a{}_i\tens f_a}
The factor $q-q^{-1}$ is an arbitrary choice of normalisation for the $e^i$,
chosen for
conventional purposes with respect to the standard examples. The normalisation
is such that the duality pairing of $B^\star,B$ is
$\ev(f_i,e^j)=\delta{}_i{}^j(q-q^{-1})^{-1}$.

For our example,  we take  $U_q(su_2)$ in more or less the standard form
generated by $q^{\pm {H\over 2}},X_\pm$ but with the opposite coproduct, and we
denote by $\widetilde{U_q(su_2)}$ its central extension by $\dila$. To present
our answer in a more familiar form, we also adjoin the square roots $q^{\pm
{H\over 4}}$ and $\dila^{\pm \h}$. We use the R-matrix form above with
ansatz\cite{FRT:lie}
\eqn{matsl2}{ \vecm^+=\pmatrix{q^{-{H\over 2}}&q^{-\h}(q-q^{-1})X_+\cr 0&
q^{H\over 2}},\quad
\vecm^-=\pmatrix{q^{H\over 2}&0\cr -q^\h(q-q^{-1})X_-&q^{-{H\over 2}}}}
(the ansatz embodies the additional relations mentioned above) and
\eqn{RR'}{ R=-\pmatrix{1&0&0&0\cr0&q^{-1}&1-q^{-2}&0 \cr
0&0&q^{-1}&0\cr0&0&0&1},\quad R'=q^2R,}
which yields the fermionic quantum planes $\C^{0|2}_q$ and
$\C^{0|2}_q{}^{\star\rm op}=\C_q^{0|2}$ as given by
\eqn{qplanes}{  \eta_+e=-q^{-1}e\eta_+,\quad e^2=\eta_+^2=0,\quad
\eta_-f=-q^{-1}f\eta_-,\quad f^2=\eta_-^2=0.}
We denote the braided (co)vector components here as
$\vece=\pmatrix{e\cr\eta_+}$ and $\vecf=(f,\eta_-)$.

\begin{propos} The double-bosonisation $\C^{0|2}_q{}^{\star\rm
op}\lbiprod\widetilde{U_q(su_2)}\rbiprod \C^{0|2}_q$ consists of $U_q(su_2)$
with the additional generators $q^{\pm {h\over 2}},\theta_\pm$ and the
relations and coproduct
\cmath{{} [q^{H\over 2},q^{h\over 2}]=0,\quad q^{H\over
2}\theta_\pm=q^{\mp\h}\theta_\pm q^{H\over 2},\quad q^{h\over 2}X_\pm=q^{\mp
\h}X_\pm q^{h\over 2}, \quad   q^{h\over 2}\theta_\pm=\pm\imath \theta_\pm
q^{h\over 2}\\
{}[\theta_+,X_-]=0,\quad [X_+,\theta_-]=0,\quad
[\theta_+,\theta_-]={q^h-q^{-h}\over q-q^{-1}},\quad
\theta_\pm^2=0,\quad [[X_\pm,\theta_\pm]_q,X_\pm]_q=0\\
\Delta q^{h\over 2}=q^{h\over 2}\tens q^{h\over 2},\quad \Delta \theta_\pm\tens
q^{-{h\over 2}}+q^{h\over 2}\tens\theta_\pm}
where $[x,y]_q=xy-qyx$ and $\imath=\sqrt{-1}$. We assume $q^2\ne\pm 1$.
\end{propos}
\proof It is instructive to see exactly  how this arises. Firstly, the quantum
group normalisation constant for (\ref{RR'}) is  $\lambda=-q^{\h}$ (to bring
$\lambda R$ into the quantum group $U_q(su_2)$ normalisation).  Next, we
compute the $\vece\vecm^+$ relations from (\ref{matdbos}), yielding
\eqn{rela}{ e q^{-{H\over 2}}=q^\h q^{-{H\over 2}}e,\quad e X_+=q^\h X_+e,\quad
\eta_+ q^{H\over 2}=q^\h q^{H\over 2}\eta_+,\quad -q^{-\h}X_+\eta_++\eta_+
X_+=q^{H\over 2}e.}
The last of these says that $e$ need not be included among the generators since
it is obtained from $\eta_+,X_+$ by a $q$-commutator. The first relation is
then implied, while the second relation remains as a `q-Serre' relation.
Similarly, the $\vece\vecm^-$ relations yield
\eqn{relb}{ X_-\eta_+=q^\h \eta_+ X_-,\quad q^\h e X_--qX_-e=\eta_+q^{H\over
2}.}
The computation for $\vecm^\pm\vecf$ is similar, with $f$ generated by
$q$-commutator of $X_-,\eta_-$.
{}From the value of $\lambda$ we obtain the cross relations
\eqn{relc}{ e\dila=-q^\h\dila e,\quad \dila f=-q^\h f\dila,\quad
\eta_\pm\dila=-q^{\pm\h}\dila\eta_\pm}
and from the form of $\vecm^\pm$ we obtain
\eqn{reld}{  [e,f]={q^{-{H\over 2}} \dila^{-1}-\dila q^{H\over 2} \over
q-q^{-1}},\quad
[\eta_+,\eta_-]={q^{H\over 2}\dila^{-1}-\dila q^{-{H\over 2}}\over
q-q^{-1}},\quad
{}[e,\eta_-]=q^{-\h}X_+\dila^{-1},\quad [\eta_+,f]=q^\h \dila X_-.}
We also obtain the coproducts
\eqn{rele}{ \Delta \eta_+=\eta_+\tens q^{H\over 2}\dila^{-1}+1\tens
\eta_+,\quad  \Delta \eta_-=\eta_-\tens 1+\dila q^{-{H\over 2}}\tens \eta_-}
and more complicated coproducts for $e,f$. The relations involving $\eta_\pm$,
along with those of $U_q(su_2)$ and $q^{h}=\dila q^{-{H\over 2}}$ provide our
new quantum group, while those with $e,f$ are redundant (they provide a useful
check). The result is a quantum group with relations $\{q^h,\eta_+\}=0$, etc.
and coproduct $\Delta\eta_+=\eta_+\tens q^{-h}+1\tens\eta_+$ etc.

Finally, if $q^{\pm {H\over 4}},\dila^{\pm \h}$ (and $\sqrt{-1}$) are also
adjoined, we can make a change of variables to
\eqn{varnew}{ \theta_+=-\eta_+q^{h\over 2},\quad \theta_-=q^{-{h\over
2}}\eta_-.}
This gives the form with a more familiar symmetrical coproduct as actually
stated. Note that the other $q$-Serre relation corresponding to
$\eta_+e=-q^{-1}e\eta_+$ and the remaining relation $e^2=0$ are implied by the
relations shown, assuming $q^2\ne -1$. Similarly for the $\eta_-,f$ sector.
\endproof

Double-bosonisation not only gives us a quantum group but constructs an
`inductive basis' and relations associated with it (such as the second of
(\ref{relb}), which are usually found using the quantum Weyl group.) The
natural basis is
\eqn{ex-basis}{ f^{\eps_1}\eta_-^{\eps_2}
\dila^{m}U_q(su_2)\eta_+^{\eps_3}e^{\eps_4},\quad \eps_i\in\{0,1\},\quad
m\in\Z,}
given a basis of $U_q(su_2)$. Similarly in the
$\theta_\pm$ version.

Also, the fermionic quantum plane algebra here is finite-dimensional and we
have a strict duality. Taking basis $\{1,e,\eta_+,\eta_+ e\}$ of $B$, the dual
basis is
\eqn{ex-dual}{\{1,(q-q^{-1})f, (q-q^{-1})\eta_-,-q(q-q^{-1})^2f\eta_-\}}
where the product is written in $B^{\star\rm op}$. Only the coefficient of
$f\eta_-$ requires some work here. It is computed from
\ceqn{ex-braid}{\ev(f\cdot_{\rm op}\eta_-,\eta_+e)=\ev(f,(\eta_+
e)\Bo)\ev(\eta_-,(\eta_+ e)\Bt)\\
 \und\Delta(\eta_+e)=\eta_+e\tens 1+1\tens\eta_+e+\eta_+\tens
e+\Psi(\eta_+\tens e),\quad \Psi(\eta_+\tens
e)=-q^{-1}e\tens\eta_+-(1-q^{-2})\eta_+\tens e}
using (\ref{matBC}). Then (\ref{dbos-R}) gives the quasitriangular structure on
the double-bosonisation as
\eqn{ex-uni}{\CR=\exp_{21}^{-1}\lambda^{-\xi\tens\xi}\CR_{su_2},\quad
\exp_{21}^{-1}=1\tens 1-(q-q^{-1})f\tens e-q^{-1}(q-q^{-1})^2 f\eta_-\tens
\eta_+e}
where we suppose that $\dila=\lambda^\xi$ so that
$\lambda^{-\xi\tens\xi}\CR_{su_2}$ is the quasitriangular structure of
$\widetilde{U_q(su_2)}$ (under which $B,B^\star$ are covariant). We used
braided-antimultiplicativity of the braided antipode to compute  $\und
S(\eta_+e)=\cdot\Psi((-\eta_+)\tens (-e))=q^{-2}\eta_+e$. Unlike the similar
`induction formula' for ${\CR}_{su_3}$ in \cite{Ma:dbos}, the braided plane
wave $\exp$ is a polynomial in our $q$-fermionic generators. On the other hand,
the dilaton contribution $\lambda^{-\xi\tens\xi}$ to the Gaussian part is more
formal. It can be treated either by adjoining $\xi$ to the algebra formally, or
by working over $\C$ with operator interpretations of all equations. Thus,
formally,
\eqn{ex-xi}{\xi={2h+H\over 2\imath\pi(\ln
q)^{-1}+1},\quad\lambda^{-\xi\tens\xi}=q^{-{(2h+H)\tens (2h+H)\over
4\imath\pi(\ln q)^{-1}+2}}}

This is the `fermionic extension' of $U_q(su_2)$ which is the cousin of
$U_q(su_3)$ as the `bosonic extension'. In the same way, all the $U_q(su_n)$
have a standard extension which yields $U_q(su_{n+1})$ (in some suitable
conventions) as in (\ref{An}), and a corresponding fermionic extension. These
are the two most natural possibilities to go from the node $U_q(su_n)$, in view
of the fact that the $su_n$ R-matrix is Hecke, i.e. has two skew-eigenvalues.
We have the same two choices for extensions of the quantum line $U_q(1)$. The
bosonic extension is $U_q(su_2)$ and the fermionic one is essentially the
non-standard quantum group studied (with quasitriangular structure) in
\cite{MaPla:non}.

\section{Anyonic and C-Statistical Doubles}

In this section we let $B$ be a $\Z/n\Z$-graded or `anyonic' braided
group\cite{Ma:any} and compute its double-bosonisation. By definition,
an anyonic braided group means a braided group which is covariant under the
anyon-generating quantum group $\Z'_{/n}$ introduced in \cite{Ma:any}. This
has a single generator $g$ and
\eqn{anygen}{ g^n=1,\quad \Delta g=g\tens g,\quad
\CR_g=n^{-1}\sum_{a,b=0}^{n-1}q^{-ab}g^a\tens g^b,}
where $q^n=1$. A covariant action means that $B$ is $\Z/n\Z$-graded, and the
corresponding (say, right) action is $b\ra g=b q^{|b|}$ for $b$ homogeneous of
degree $|b|$.
The induced braiding is $\Psi(b\tens c)=c\tens b q^{|c||b|}$ for homogeneous
$b,c\in B$, and $B$ is required to be a braided group with this braiding. We
also need another braided group $B^\star$ dual to $B$ in the sense of a map
$\ev$ relating the (opposite) product of one  to the coproduct of the other,
and covariant in the sense $\ev(c,b)=0$ unless $|c|+|b|=0$ for all homogeneous
$c\in B^\star$ and $b\in B$. Putting these formulae into Theorem~2.1, we obtain
at once:

\begin{propos} For  $B$ anyonic with dual $B^\star$, the double-bosonisation
$B^{\star\rm op}\lbiprod \Z'_{/n}\rbiprod B$ is a quantum group containing
$B^{\star\rm op},\Z'_{/n},B$ with cross relations and coproduct
\cmath{ bg=gb q^{|b|},\quad c g= g c q^{|c|},\quad b\Bo g^{|c\Bt|} c\Bo
\ev(c\Bt,b\Bt)=\ev(c\Bo,b\Bo)c\Bt g^{|b\Bo|}b\Bt\\
\Delta b=b\Bo\tens g^{|b\Bo|}b\Bt,\quad \Delta c=g^{|c\Bt|} c\Bo\tens c\Bt}
for all $b\in B$, $c\in B^\star$.\end{propos}

These computations are similar to the derivation of the formulae for single
bosonisation in \cite{Ma:book}. In the case of strict duality, we also have a
quasitriangular structure $\CR=\exp_{21}^{-1}\CR_g$ from (\ref{dbos-R}). On the
other hand, when $n=2$ any anyonic braided group means a super-quantum group.
In this case, we find that we can recognise the double-bosonisation as the
usual bosonisation of a suitable super-double. To this end, we first define the
anyonic-opposite product on $B^\star$, which algebra we denote
$B^{\star\und{\rm op}}$, by:
\eqn{diam}{ c\cdot_{\und{\rm op}} d=dc q^{-|c||d|}.}

\begin{propos} When $n=2$ and $q=-1$, we have $B^{\star\rm op}\lbiprod
\Z'_{/n}\rbiprod B\isom \Z'_{/n}\rbiprod D(B)$, where $D(B)$ is generated by
$B^{\star\und{\rm op}},B$ as subalgebras, with the cross relations
\[ b\Bo c\Bo
\ev(c\Bt,b\Bt)(-1)^{|b\Bo||c|}=(-1)^{|c\Bo||c\Bt|}\ev(c\Bo,b\Bo)c\Bt b\Bt\]
for all $b\in B$, $c\in B^\star$. The coproducts are those of $B^\star,B$.
\end{propos}
\proof The novel part here is not the definition of $D(B)$, which is analogous
to Drinfeld's quantum double (up to signs), but the isomorphism
$\theta:\Z'_{/n}\rbiprod D(B)\to B^{\star\rm op}\lbiprod \Z'_{/n}\rbiprod B$,
which we define as
\[ \theta(g)=g,\quad \theta(b)=b,\quad \theta(c)=g^{|c|}c,\quad
\forall b\in B,\ c\in B^{\star\und{\rm op}}\]
This is clearly an isomorphism when restricted to the subalgebras
$B^{\star\und{\rm op}},\Z'_{/n},B$, and respects the $g-b$ and $g-c$ cross
relations in $\Z'_{/n}\rbiprod D(B)$ since these have the same form as in
Proposition~4.1. Moreover,
applying $\theta$ to the left hand side of the $D(B)$ relation gives
\align{&&(-1)^{|b\Bo||c|}b\Bo g^{|c\Bo|}c\Bo\ev(c\Bt,b\Bt)=g^{|c|}b\Bo
g^{|c\Bt|}c\Bo\ev(c\Bt,b\Bt)=g^{|c|} \ev(c\Bo,b\Bo)c\Bt g^{|b\Bo|}b\Bt\\
&&\qquad=g^{|c|} \ev(c\Bo,b\Bo)c\Bt
g^{|c\Bo|}b\Bt=(-1)^{|c\Bo||c\Bt|}\ev(c\Bo,b\Bo)g^{|c\Bt|}c\Bt b\Bt,}
which is $\theta$ applied to the right hand side. The coproducts also map over:
the coproduct of $\Z'_{/n}\rbiprod D(B)$ is the same as in Proposition~4.1 on
$B,\Z'_{/n}$, and is $\Delta c=c\Bo\tens g^{|c\Bo|}c\Bt$ on
$B^{\star\und{\rm op}}$. We have
$(\theta\tens\theta)\circ\Delta c=g^{|c\Bo|}c\Bo\tens
g^{|c|}c\Bt=g^{|c|}g^{|c\Bt|}c\Bo\tens g^{|c|}c\Bt=\Delta\circ\theta(c)$.
\endproof

When $n>2$, we still obtain an ordinary quantum
group by Proposition~4.1 but it is not (as far as I know) the single
bosonisation of anything. This is because the braided quantum double
construction does not really work (it gets tangled up) and, indeed, the
double-bosonisation is probably the closest we can come to it. This was one of
the motivations in \cite{Ma:dbos}.

Finally, we are not limited to one-dimensional gradings. In \cite{Ma:csta} we
introduced quasitriangular Hopf algebras $U_q(\beta)$ (say) associated to any
bilinear form $\beta_{ij}$. The generators are commuting variables $\{\xi_i\}$
with
structure
\eqn{csta}{ \Delta \xi_i=\xi_i\tens 1+1\tens \xi_i,\quad \CR_\xi=q^{\sum
\beta_{ij}\xi_i\tens \xi_j},}
where $q$ is general. The action on a $\Z^n$-graded algebra is $b\ra \xi_i=
b|b|_i$ where $|\ |_i$ is the $i$-component of the degree. Then $\Psi(b\tens
c)=c\tens b q^{\sum \beta_{ij}|c|_i|b|_j}$ is the braiding. Braided groups in
this setting are called $\C$-statistical\cite{Ma:csta}.

An example of $B,B^\star$ in this setting is provided\cite{Ma:csta} by the free
algebras
$B=\C\<e^i\>$, $B^\star=\C\<f_i\>$, with grading  and resulting braiding and
pairing
\eqn{ex-free}{|e^i|_j=\delta^i{}_j,\quad \Psi(e^i\tens e^j)=e^j\tens e^i
q^{\beta_{ji}},\quad |f_i|_j=-\delta_{ij},\quad \Psi(f_i\tens f_j)=f_j\tens
f_iq^{\beta_{ji}},\quad \ev(f_i,e^j)=\delta_i{}^j(q^i-q_i^{-1})^{-1}.}
This is simply the R-matrix setting (\ref{matBC}) with
$R^i{}_j{}^k{}_l=\delta^i{}_j\delta{}^k{}_l q^{\beta_{ki}}$ and $R'=P$, the
permutation matrix. (More generally, one can have {\em free braided groups}
associated to any $R$ and $R'=P$, as studied extensively in \cite{Ma:fre}.)
The additional parameters denoted $q^i-q_i^{-1}$ reflect some freedom in the
normalisation of the $e^i$. From Theorem~2.1, we have immediately:

\begin{propos} The double-bosonisation $\C\<\vecf\>\lbiprod
U_q(\beta)\rbiprod\C\<\vece\>$ has mutually non-commuting $\{e^i\}$, mutually
non-commuting $\{f_j\}$, mutually commuting $\{\xi_i\}$ and cross relations and
coproduct obeying
\cmath{ e^iq^{H^j}=q^{H^j} e^i q^{\beta_{ji}},\quad f_iq^{H_j}=q^{H_j} f_i
q^{-\beta_{ij}},\quad [e^i,f_j]=\delta^i{}_j{q^{H^i}-q^{-{H_i}}\over
q_i-q_i^{-1}}\\
\Delta e^i=e^i\tens q^{H^i}+1\tens e^i,\quad \Delta f_i=f_i\tens
1+q^{-{H_i}}\tens f_i;\quad H^i\equiv \sum\beta_{ij}\xi_j,\quad H_i\equiv \sum
\xi_j\beta_{ji}.}
\end{propos}

The computation is similar to the single bosonisation of (\ref{ex-free}) in
\cite{Ma:csta}. We assume for simplicity here that $\beta$ is invertible --
otherwise we have to work with the $\xi$ variables, with relations
$[e^i,\xi_j]=\delta^i{}_je^i$ and $[f_i,\xi_j]=-\delta_{ij}f_i$. Also,
we have presented the $H^i,H_i$ symmetrically, but either set will suffice.

The pairing $\ev$ between the free braided groups (\ref{ex-free}) is typically
degenerate, but we can always quotient out by the kernels of the pairing on
each side. When $\beta_{ij}$ is the symmetric bilinear form associated to a
Cartan matrix and $q^i=q_i=q^{\beta_{ii}\over 2}$, Lusztig\cite{Lus} has
effectively computed the kernels and found that they are generated by the
$q$-Serre relations. Hence, quotienting by these gives $\C$-statistical braided
groups $U_q(n_\pm)$ and
\eqn{Lus}{ U_q(n_-)\lbiprod U_q(\beta)\rbiprod U_q(n_+)=U_q(\cg)}
recovers Lusztig's construction of $U_q(\cg)$ in our braided approach as a
double-bosonisation. We have glossed over a lot of technicalities here (we need
to work with $q^{H_i}$ as generators and avoid the power series inherent in
$\CR_\xi$, by working with weak quasitriangular structure maps); see
\cite{Ma:dbos} for the full details.

On the other hand, as explained in \cite{Ma:dbos}, we are not limited to
$\beta_{ij}$ associated to a Cartan matrix. When we take some other bilinear
form (typically integer-valued, for an algebraic answer) we have some other
kernel of $\ev$ in (\ref{ex-free}). The extension of $\ev$ to products is via
braided differentiation, so the kernels have a `braided geometrical'
interpretation\cite{Ma:fre}\cite{Ma:book}. Quotienting out by the kernels, we
obtain  non-degenerately paired braided groups $B,B^\star$. Their
double-bosonisation gives new quantum groups complete with quasitriangular
structure $\exp_{21}^{-1}\CR_\xi$, where $\exp\in B\tens B^\star$ is the
braided exponential or coevaluation for the pairing $\ev$, generally existing
as a formal power series. Some of these quantum groups will be computed
elsewhere.

%\bibliographystyle{unsrt}
%\bibliography{biblio}

\end{document}